\documentclass[12pt]{article}

\usepackage{amsmath,amssymb,latexsym}

\newtheorem{defn}{Definition}[section]
\newtheorem{lemma}[defn]{Lemma}
\newtheorem{ex}[defn]{Example}
\newtheorem{thm}[defn]{Theorem}
\newtheorem{prop}[defn]{Proposition}

\newtheorem{rem}[defn]{Remark}

\newcommand{\h}{{H}}

\newcommand{\mn}{\mathbb N}
\newcommand{\mr}{\mathbb R}

\newcommand{\mprop}{\mathcal{P}}
\newcommand{\sspace}{{\mathfrak{s}}}
\newcommand{\bs}{{\bf s}}

\newcommand{\spaceh}{{\mathfrak{H}}}
\newcommand{\nulel}{{\bf 0}}

\newcommand{\seqgr[1]}{(#1_n)_{n=1}^\infty }

\newcommand{\sech[1]}{{\cap_{k\in\mn_0}{#1}_k}}

\newcommand{\snorm[1]}{{\|\hspace{-.01in}|{#1}\|\hspace{-.015in}|}}

\def\be{\begin{equation}}
\def\ee{\end{equation}}
\def\newin {\kern-0.22em\in\kern-0.15em}
\def\newsubset {\kern-0.2em\subset\kern-0.2em}

\def\<{\langle}
\def\>{\rangle}

\begin{document}

\title{Frame expansions of test functions,  tempered distributions, and  ultradistributions}

\author{\vspace{.1in} Stevan Pilipovi\'c$^{ \, a}$  and Diana T. Stoeva$^{\, b}$ \\
{}
 \vspace{-.07in} {\small 
$^{a}$ Department of Mathematics and Informatics, University of Novi Sad,} \\
\vspace{-.07in} 
{\small   Trg D. Obradovi\'ca 4, 21000 Novi Sad, Serbia, pilipovic@im.ns.ac.yu} \\
\vspace{-.07in} 
 {\small  $^{b}$ Acoustics Research Institute, Austrian Academy of Sciences,}
\\
\vspace{-.07in} 
{\small  Wohllebengasse 12-14, Vienna 1040, Austria,  dstoeva@kfs.oeaw.ac.at
}
}

\maketitle

\abstract{
The paper is devoted to frame expansions in Fr\'echet spaces. First we  review  some results which concern series expansions in general Fr\'echet spaces via Fr\'echet and General Fr\'echet frames. Then we present some new results on series expansions of  tempered distributions and ultradistributions, and the corresponding test functions, via localized frames and coefficients in appropriate sequence spaces. 
}

\vspace{.1in} 

\vspace{.05in} 

\section{Introduction}

In this paper we present results devoted to frame expansions in Fr\'echet spaces. First we consider the general case, expansions via Fr\'echet frames and appropriate dual sequences in general, and then we aim expansions of generalized functions via a proper class of frames. 
As the Hermite expansions
are the basic ones for tempered distributions and ultradistributions, a suitable localization of a frame
with respect to the Hermite basis enables us to
extend the consideration of generalized functions using appropriate class of frames instead of the Hermite basis. Frames were introduced in Hilbert spaces \cite{DSframe}. They generalize the concept of an orthonormal basis allowing even redundancy, but still provide series expansions of all the elements of the space. Frames were extended  to Banach spaces (atomic decompositions and Banach frames \cite{FGat1,FGat2,Gbanach}, p-frames \cite{AST}, $X_d$-frames \cite{CCS})  and furthermore to Fr\'echet spaces (pre-Fr\'echet, Fr\'echet, and General Frechet frames, \cite{pst,ps2,ps4}). While Hilbert frames always guarantee series expansions in Hilbert spaces, this is not always the case with Banach and Fr\'echet frames.  
In this paper we review some results related to sufficient conditions for series expansions in general Fr\'echet spaces, as well as present new results  devoted to series expansions in certain spaces of test functions and their duals via appropriate frames. 
The paper is organized as follows. Section \ref{sec-Prel} contains the main definitions, notation, and basic needed facts. In Section \ref{sec3} we review some results from \cite{ps2,ps4} which are related to series expansions in general Fr\'echet spaces via Fr\'echet and General Fr\'echet frames. Section \ref{sec4} is devoted to  new results which concern expansions in certain spaces of test functions, tempered distributions, and ultradistributions, via localized frames; the statements are given without proofs and 
the proofs are subject of a further extended paper \cite{ps5}.

\section{Preliminaries}\label{sec-Prel}

\vspace{.1in}
Throughout the paper, 
$(\h, \<\cdot,\cdot\>)$  denotes a separable Hilbert space. We  consider countable sequences and for convenience of the writing we index them by the set $\mn$. 
A sequence $\seqgr[g]$ with elements from $\h$ is called: a
 \emph{frame for $\h$} 
 if there exist positive constants $A$ and $B$ (called {\it frame bounds}) so that 
$A\|f\|^2 \leq \sum_{n=1}^\infty |\<f,g_n\>|^2\leq B\|f\|^2$ for every $f\in\h$ \cite{DSframe}; a \emph{Riesz basis for $\h$} 
 if its elements are the images of the elements of an orthonormal basis under a bounded bijective operator on $\h$ \cite{Bari51}. 
Let us recall some needed basic facts from frame theory  (see e.g \cite{Cbook}). 
Let $G=\seqgr[g]$ be a frame for $\h$. 
Then there exists a frame $\seqgr[f]$ for $\h$ so that 
$f=\sum_{n=1}^\infty \<f,f_n\>g_n=\sum_{n=1}^\infty \<f,g_n\>f_n, f\in\h.$  
Such $\seqgr[f]$ is called a {\it dual frame of $\seqgr[g]$}. 
The {\it analysis operator} $U_G$, given by $U_Gf=(\<f,g_n\>)_{n=1}^\infty$, is bounded from $\h$ into $\ell^2$, the {\it synthesis operator} $T_G$ defined on the finite sequences 
by $T_G(c) = \sum_n c_n e_n$ extends to a bounded operator from $\ell^2$ into $\h$, the {\it frame operator} $S_G=T_G U_G$ 
is a bounded bijection of $\h$ onto $\h$, and the series in $S_Gf=\sum_{n=1}^\infty \<f, g_n\> g_n$ converges unconditionally. 
The sequence $(S_G^{-1}g_n)_{n=1}^\infty$ is a dual frame of $\seqgr[g]$,  called the {\it canonical dual of $\seqgr[g]$},   and it will be denoted by $(\widetilde{g_n})_{n=1}^\infty$.
 When $\seqgr[g]$ is a Riesz basis of $\h$ (and thus a frame for $\h$), then 
 $(\widetilde{g_n})_{n=1}^\infty$ is the only dual frame of $\seqgr[g]$. 
 Frames which are not Riesz bases have other dual frames in addition to the canonical dual. 
\sloppy
Recall also the localization-notions introduced in \cite{G2004}. Given a Riesz basis $\seqgr[g]$ for  $\h$, a frame $E=\seqgr[e]$ for $\h$ is called: 
 {\it polynomially localized with respect to $\seqgr[g]$ with decay $s>0$} 
if there is a constant $C_s>0$ so that 
$
 \max\{ |\<e_m, g_n\>|, |\<e_m, \widetilde{g_n}\>|\} \leq 
 C_s(1+|m-n|)^{-s}, \ m,n\in\mn;$
 {\it exponentially localized with respect to $\seqgr[g]$} if for some $s>0$ there is a constant $C_s>0$ so that
$
 \max\{ |\<e_m, g_n\>|, |\<e_m, \widetilde{g_n}\>|\} \leq 
 C_s   \mathrm{e}^{-s |m-n|}, \ m,n\in\mn.
$ Notice that in the literature there exist other ways to define 
 localization of frames \cite{BCHL1, BCHL2, FG, Futamura}, but for the purpoces of the current paper it is relevant and enough to use the localization concepts according  to \cite{G2004}.

\vspace{.1in}
Next, $(X, \|\cdot\|)$ denotes a Banach space and 
 $( \Theta, \snorm[\cdot]) $ denotes a Banach sequence space. 
A Banach sequence space is a {\it $BK$-space} if the coordinate functionals are continuous.
If the canonical vectors form a Schauder basis for $\Theta$, then $\Theta$ is called a {\it $CB$-space} and it is clearly a $BK$-space. 
 Given a $BK$-space $\Theta$ and a Riesz basis $G=(g_n)_{n=1}^\infty$ for $\h$,  one associates to $\Theta$ the following Banach space 
$$\spaceh^\Theta_G:=\{f\in\h \ : \ f=\sum_{n=1}^\infty c_n g_n 
\mbox{ with $(c_n)_{n=1}^\infty\in\Theta$} \} \, \mbox{ normed by } \, \|f\|_{\spaceh^\Theta_G}:=\snorm[ (c_n)_{n=1}^\infty]_{\Theta}.$$

 Further, we consider Fr\'echet spaces which are projective limits of Banach spaces. 
Let $\{Y_k, | \cdot |_k\}_{k\in\mn_0}$ be a sequence of separable
Banach spaces such that \be \label{fx1} \{\nulel\} \neq
\sech[Y]\subseteq \ldots \subseteq Y_2 \subseteq Y_1 \subseteq Y_0
\ee \be  \label{fx2} |\cdot|_0\leq | \cdot |_1\leq | \cdot |_2\leq
\ldots \ee \be \label{fx3} Y_F :=\sech[Y] \;\; \mbox{is dense in}
\;\; Y_k, \;\; k\in\mn_0. \ee 
Under the conditions (\ref{fx1})-(\ref{fx3}), $Y_F$ is a Fr\'echet space with the sequence of norms $ |\cdot |_s, $ $ s\in\mn_0$, and it is called the {\it projective limit} of $Y_k$, $k\in\mn_0$. 
 We will use such type of sequences in two cases:

1. $Y_k=X_k$ with norm $\|\cdot\|_k, k\in\mn_0;$

2. $Y_k=\Theta_k$ with norm $\snorm[\cdot]_k, k\in\mn_0$.

\vspace{.1in}
 We use the therm \emph{operator} for a linear mapping.  Given sequences of Banach spaces, $\{X_k\}_{k\in\mn_0}$ and $\{\Theta_k\}_{k\in\mn_0}$,
which satisfy (\ref{fx1})-(\ref{fx3}), an operator $G: \Theta_F \to X_F$
is called {\it $F$-bounded} if for every 
$k\newin\mn_0$, 
there exists a constant $C_k>0$ such that $\|G\seqgr[c]\|_k\leq
C_k \snorm[\{c_n\}_{n=1}^\infty]_k$ for all $\seqgr[c]\newin
\Theta_F$. Now we recall the definitions of Fr\'echet and General  Fr\'echet frames.

\begin{defn}\label{fframe}  {\rm  \cite{ps4}} 
Let $\{X_k, \|\cdot\|_k\}_{k\in\mn_0}$ be a sequence of Banach
spaces satisfying (\ref{fx1})-(\ref{fx3}) and let
$\{\Theta_k, \snorm[\cdot]_k\}_{k\in\mn_0}$ be a sequence of
$BK$-spaces satisfying (\ref{fx1})-(\ref{fx3}). A sequence
$\seqgr[g]$ with elements from ${X_F^*}$ is called: a {\it General  pre-Fr\'echet frame} (in short, {\it General pre-$F$-frame})  for $X_F$ with respect to $\Theta_F$ if 
there exist sequences $\{\widetilde{s}_k\}_{k\in\mn_0}\subseteq \mn_0$, $\{s_k\}_{k\in\mn_0}\subseteq \mn_0$, which increase to $\infty$ with the property $s_k\leq \widetilde{s}_k$, $k\in\mn_0$, and there
exist constants $0<A_k\leq B_k<\infty$, $k\in\mn_0$, satisfying 
\begin{equation}\label{fframestar}
(g_n(f))_{n=1}^\infty\in\Theta_F, \ f\newin X_F,
\end{equation}
\begin{equation}\label{fframetwostar}
A_k \|f\|_{s_k}\leq \snorm[\{g_n(f)\}_{n=1}^\infty]_k\leq
B_k\|f\|_{\widetilde{s}_k}, \ f\newin X_F, k\in\mn_0;
\end{equation}
 a {\it General  Fr\'echet frame} (in short, {\it General $F$-frame})  for $X_F$ with respect to $\Theta_F$ if it is a General pre-$F$-frame for $X_F$ with respect to $\Theta_F$ 
and there exists a continuous operator 
$V:\Theta_F\rightarrow X_F$ so that $V(g_n(f))_{n=1}^\infty=f$
for every $f\newin X_F$.
\end{defn}

Let $s_k=\widetilde{s}_k=k$, $k\in\mn_0$. In this case 
the above definition of a General pre-$F$-frame reduces to the definition of a  {\it pre-Fr\'echet frame} (in short {\it pre-$F$-frame}) \cite{pst},
and if in addition the continuity of $V$ is replaced by the stronger condition of $F$-boundedness of $V$, then one comes to the concept of a  {\it Fr\'echet frame} (in short {\it $F$-frame}) 
\cite{ps2}. 
 In the particular case when $X_k=X$, and $\Theta_k=\Theta$,  $k\in\mn_0$:   
 an $F$-frame (resp. pre-$F$-frame)  for $X_F$ with respect to $\Theta_F$ is actulally a {\it Banach frame for $X$ with respect to $\Theta$} (resp. $\Theta$-frame for $X$) as introduced in \cite{Gbanach}  (resp. \cite{CCS});
 when (\ref{fframestar}) and the upper inequality of (\ref{fframetwostar}) hold, one comes to the definition of a {\it $\Theta$-Bessel sequence for $X$}.

When $(g_n)_{n=1}^\infty$ is a pre-F-frame for $X_F$ with respect to $\Theta_F$, then for any $n\in\mn$ and any $k\in\mn_0$, $g_n$ can be extended in a unique way to a continues operator on $X_k$ and this extension will be denoted by $g_n^k$.

\vspace{.1in} Recall that a positive
continuous 
function $\mu$ on $\mr$ is called: a  {\it $k$-moderate weight} if $k\geq 0$ and there exists a  constant $C>0$  so that
$ \mu(t+x)\leq C (1+|t|)^k \mu(x), \ t,x\in\mr;$
a  {\it sub-exponential}  {weight}, if there exist  constants  $C>0,\gamma >0$ and $\beta\in (0,1)$  so that
$ \mu(t+x)\leq C e^{\gamma |t|^\beta}\mu(x), \  t,x\in\mr.$ 
Taking $\mu_k(x)=(1+|x|)^k$ (resp. $\beta\in(0,1)$ and $\mu_k(x)=e^{k|x|^\beta}$),   $ x\in\mr$,  $k\in\mn_0$, the spaces $\Theta_k:=\ell^2_{\mu_k}$, $k\in\mn_0$, satisfy (\ref{fx1})-(\ref{fx3}) and their projective limit $\cap_k \Theta_k$ is the so called  {\it space of rapidly decreasing sequences}
$\bs$ (resp. {\it space of sub-exponentially decreasing sequences} $\sspace^{\beta}$).
 Further, we will use the following statement:

 \begin{lemma} \label{lemrb} 
Let $G=\seqgr[g]$ be a Riesz basis for $\h$. 
For $k\in\mn_0$, 
let  $\mu_k$ be a  $k$-moderate weight  
so that   
$1=\mu_0(x)\leq \mu_1(x) \leq \mu_2(x) \leq ...$, for every $x\in\mr$. 
Then   $\{\Theta_k\}_{k\in\mn_0}:=\{\ell^2_{\mu_k}\}_{k\in\mn_0}$ is a sequence of $CB$-spaces satisfying (\ref{fx1})-(\ref{fx3}),  
 the spaces 
  $X_k:=\spaceh^{\Theta_k}_G$, $k\in\mn_0$, satisfy (\ref{fx1})-(\ref{fx3}),  and
  $g_n\in X_F$ for every $n\in\mn$. 
 The conclusions also hold if the assumptions ``$\mu_k$ -  $k$-moderate weight'' are replaced by ``$\mu_k$ - sub-exponential weight''. 
\end{lemma}

Test function spaces and their duals under consideration in the paper are:\\ 
${\mathcal S} (\mr) := \{ f\in C^\infty(\mr) \ :  |f|_{k}:=\sup_{x\in\mr} \sup_{m\leq k} |f^{(m)}(x)|\, (1+|x|^2)^{k/2}<\infty, \  \forall k\in\mn_0 \}$,
and its dual $ {\mathcal S}' (\mr)\subset C^\infty(\mr)$, the space of tempered distributions;\\
$\Sigma^\alpha:=\{\phi\in C^{\infty}(\mr)\ :  |\phi|_{h,\alpha}: = 
\sup_{n\in\mn_0, x\in\mr} \frac{h^{n}e^{m|x|^{1/\alpha}}|\phi^{(n)}(x)|}{ n!^\alpha}<\infty, \ \forall h>0\},$
and its dual $ ({\Sigma^\alpha}(\mr))'$, $\alpha>1/2,$ the space of Beurling  tempered ultradistributions, cf.  \cite{sp,toft,toft2}.

In the sequel,  $(h_n)_{n=1}^\infty$ is the Hermite orthonormal basis $({\bf h}_n)_{n=0}^\infty$ of $L^2(\mr)$, re-indexed from $1$ to $\infty$, i.e., 
 $h_{n+1}(t)={\bf h}_{n}(t) 
= (2^{(n)} n!\, \sqrt{\pi})^{-1/2}  \, (-1)^{n} e^{t^2/2} \frac{d^{n}}{dt^{n}} (e^{-t^2})$, $n\in\mn_0$.  
Recall that $h_n\in\mathcal{S}$ and $h_n\in\Sigma^\alpha$, $\alpha>1/2$, $n\in\mn$.
Moreover, we know
\cite{Simon} the following:

- If $f\in\mathcal{S}$, then $(\<f, h_n\>)_{n=1}^\infty\in\bs$; conversely, if $(a_n)_{n=1}^\infty\in\bs$, then $\sum_{n=1}^\infty a_n h_n$ converges in $\mathcal{S}$ to  $f=\sum_{n=1}^\infty\<f, h_n\>h_n, \;(\<f, h_n\>)_{n=1}^\infty=(a_n)_{n=1}^\infty$. 

-  If  $F\in \mathcal{S}'$, then 
$(b_n)_{n=1}^\infty:=(F(h_n))_{n=1}^\infty\in\bs'$ and $F(f)=\sum_{n=1}^\infty \<f, h_n\> b_n$,  $f\in  \mathcal{S}$; 
conversely, if $(b_n)_{n=1}^\infty\in\bs'$, then the mapping $F: f\to \sum_{n=1}^\infty \<f, h_n\> b_n$ 
 is well defined on $\mathcal{S}$, it determines $F$ as an element of $\mathcal{S}'$ 
and $(F(h_n))_{n=1}^\infty=(b_n)_{n=1}^\infty$.

The above two statements also hold when $\mathcal{S}$, $\mathcal{S}'$, $\bs$, and $\bs'$ are replaced by 
$\Sigma^\alpha$, $(\Sigma^\alpha)'$, $\sspace^{1/(2\alpha)}$, and $(\sspace^{1/(2\alpha)})'$ with $\alpha>1/2$, 
respectively  \cite{sp,toft,toft2}.

\section{Frame expansions in Fr\'echet spaces}\label{sec3}

In this section we recall some general statements about sufficient conditions for series expansions in Fr\'echet spaces via Fr\'echet and General Fr\'echet frames, and appropriate dual sequences. We start with the case of  Fr\'echet frames.

\begin{prop} \label{diff} {\rm \cite{ps2}} Let
$\seqgr[g]$ be an $F$-frame for $X_F$ with respect to $\Theta_F$.
Then the following holds.
\begin{itemize}
\item[{\rm(a)}] For every $k\newin\mn_0$, the sequence $\{g_i^k\}_{i=1}^\infty$
is a Banach frame for $X_s$ with respect to $\Theta_k$.
\item[{\rm(b)}] If \,$\Theta_k$, $k\newin\mn_0$,
are  $CB$-spaces, then there exists  $\seqgr[f]\newin (X_F)^\mn$,
which is $\Theta_k^*$-Bessel sequence for $X_k^*$ for every $k\in\mn_0$ and such that
\begin{eqnarray}
f&=&\sum_{i=1}^{\infty} g_i(f) f_i,  \   f\newin X_F, \
\mbox{(in $X_F$),} \label{frepr} \\
g&=&\sum_{i=1}^{\infty} g(f_i) g_i,  \  g\in X_F^*, \
\mbox{(in $X_F^*$),} \label{frepr2}\\
f&=&\sum_{i=1}^{\infty} g_i^k(f) f_i, \  f\newin X_k, \
 k\newin\mn_0.\label{fsrepr}
\end{eqnarray}
\item[{\rm(c)}] If \,$\Theta_k$ and  $\Theta^*_k$, $k\newin\mn_0$, are  $CB$-spaces,
then there exists $\seqgr[f]\newin (X_F)^\mn$, which is a
$\Theta_k^*$-frame for $X_k^*$ for every $k\in\mn_0$ and such that
(\ref{frepr})-(\ref{fsrepr}) hold, and moreover,
\begin{equation} \label{xs}
g=\sum_{i=1}^{\infty} g(f_i) g_i^k, \ g\in X_k^*,\,
 k\in\mn_0. \ \ \ \ \ \ \ \ \ \ \ \ \ \ \ \ \ \ \ \ \ \ \ \ \ \ \
\end{equation}
\item[{\rm(d)}]
If \,$\Theta_k$, $k\newin\mn_0$, are reflexive $CB$-spaces, then
there exists $\seqgr[f]\newin (X_F)^\mn$, which is a Banach frame
for $X_k^*$ with respect to $\Theta_k^*$ for every $k\in\mn_0$ such that
(\ref{frepr})-(\ref{xs}) hold.
\end{itemize}
\end{prop}

As one can see in Proposition \ref{diff},  
the F-boundedness property of $V$ 
 leads to series expansions in all the spaces $X_k$, $k\in\mn_0$. 
 Now we continue with the case of General Fr\'echet frames and show that in this case the continuity property of $V$ is enough to imply the existence of a subsequence 
 $\{X_{\widetilde{w}_{j}}\}_{j=0}^\infty$ of the given sequence $\{X_{\widetilde{s}_k}\}_{k=0}^\infty$ according to Definition \ref{fframe}, so that one has series expansions in $X_{\widetilde{w}_{j}}$, $j\in\mn_0$, with convergence in appropriate norms. 
 
\begin{thm} \label{th2} {\rm \cite{ps4}}
 Let $\seqgr[g]$ be a General $F$-frame for $X_F$ with respect to $\Theta_F$ and let $\Theta_{k}$, $k\newin\mn_0$, be  $CB$-spaces. Then there exist  sequences $\{w_j\}_{j\in\mn_0}$, $\{r_j\}_{j\in\mn_0}$, and $\{\widetilde{w}_j\}_{j\in\mn_0}$, which increase to $\infty$ and there exist constants $\widetilde{A}_j, \widetilde{B}_j$, $j\in\mn_0$, such that
 for every $j\in\mn_0$, 
\begin{equation*}
\widetilde{A}_j \|f\|_{w_j}\leq \snorm[\{g_i(f)\}_{i=1}^\infty]_{r_j}\leq
\widetilde{B}_j\|f\|_{\widetilde{w}_j}, \ \forall f\newin X_F.
\end{equation*} 
  Moreover, there exists a sequence  $\seqgr[f]\newin (X_F)^\mn$ such that for every $j\in\mn_0$, $\seqgr[f]$
 is a $\Theta_{r_j}^*$-Bessel sequence for $X_{w_j}^*$ and
 \begin{eqnarray*}  
f&=&\sum_{i=1}^{\infty} g_i^{\widetilde{w}_j}(f) f_i \ \mbox{in $\|.\|_{w_j}$-norm}, \  f\newin X_{\widetilde{w}_{j}}.
\end{eqnarray*}

\end{thm}

\section{Expansions of tempered distributions and ultradistributions by localized frames} \label{sec4}

In this section, we aim 
 expansions in a Fr\'echet space and its dual by 
 localized frames and coefficients in a corresponding Fr\'echet sequence space.  
 First we present a general result based on frames localized with respect to a Riesz basis and providing frame expansions in a corresponding Fr\'echet space,
 and then we apply it  
 to obtain frame expansions 
in spaces of tempered distributions and ultradistributions, as well as in the corresponding test function spaces.  
To clarify some notation, when we take a frame element $e_n\in X_F(\subset \h\subset X_F^*)$ and consider it as a functional in $X_F^*$, then we denote it by bold-style ${\bf e_n}$.

Let us start with the general theorem about expansions in Fr\'echet spaces via localized frames.  

\begin{thm} \label{proprb1}
Let the assumptions and notation of Lemma \ref{lemrb} hold. 
Assume that $E=(e_n)_{n=1}^\infty$  is a sequence with elements from $X_F$ which is a frame for $\h$ and polynomially localized with respect to $G$ with decay $s$ for every $s\in\mn$  (resp. exponentially  localized with respect to $G$).   
Then the following statements hold.
 \begin{itemize}
 \item[{\rm (a)}]  $\widetilde{e_n}\in X_F$ for every $n\in\mn$.
\item[{\rm (b)}]    
The analysis operator $U_E$ is $F$-bounded from $X_F$ into $\Theta_F$, the synthesis operator $T_E$ is $F$-bounded from $\Theta_F$ into $X_F$, the frame operator $S_E$ is $F$-bounded and bijective from $X_F$ onto $X_F$ with unconditional convergence of the series in $S_Ef=\sum_{n=1}^\infty \<f, e_n\> e_n$, and $\widetilde{e_n}\in X_F$, $n\in\mn$. 
 
   \item[{\rm (c)}] For every $f\in  X_F$, 
  $$\mbox{$f=\sum_{n=1}^\infty \<f, \widetilde{e_n} \> e_n = \sum_{n=1}^\infty \<f, e_n\>\widetilde{e_n}$ (with convergence in $X_F$)} $$
with $(\<f,  \widetilde{e_n} \>)_{n=1}^\infty\in\Theta_F$ and $(\<f,  e_n \>)_{n=1}^\infty\in\Theta_F$.

\item[{\rm (d)}] If $X_F$ and $\Theta_F$ have the following property with respect to $(g_n)_{n=1}^\infty$:
{\rm $$ \mbox{$\mprop_{(g_n)}$: For  $f\in \h$, one has $f\in  X_F$ if and only if $(\<f,  g_n \>)_{n=1}^\infty\in\Theta_F$.}$$}
 then they also have the properties $\mprop_{(e_n)}$ and $\mprop_{(\widetilde{e_n})}$.

 \item[{\rm (e)}] 
Both  $({\bf e_n})_{n=1}^\infty$ and $(\widetilde{\bf e_n})_{n=1}^\infty$   form Fr\'echet frames for $X_F$ with respect to $\Theta_F$.

 \item[{\rm (f)}]  
  For every $g\in  X_F^*$, 
\begin{equation}\label{reprg}
\mbox{$g    = \sum_{n=1}^\infty g(e_n) \, \widetilde{{\bf e_n}}
    =\sum_{n=1}^\infty g(\widetilde{e_n}) \, {\bf e_n}$ (with convergence in $X_F^*$)} 
    \end{equation}
with $(g(e_n) )_{n=1}^\infty\in\Theta_F^*$ and $(g(\widetilde{e_n}) )_{n=1}^\infty\in\Theta_F^*$.

 \item[{\rm (g)}]    If $(a_n)_{n=1}^\infty \in \Theta_F^*$, then $\sum_{n=1}^\infty a_n {\bf e}_n$ (resp. $\sum_{n=1}^\infty a_n \widetilde{\bf e_n}$) converges in $X_F^*$, i.e., the mapping $f\mapsto \sum_{n=1}^\infty \<f, e_n\> a_n$  (resp. $f\mapsto \sum_{n=1}^\infty \<f, \widetilde{e_n}\> a_n$) determines a continuous linear functional on $X_F$. 

\end{itemize} 

\end{thm}

\noindent 
{\bf Remark:} Note that in the setting of the above theorem, when $G$ is an orthonormal basis of $\h$ or more generally, when $G$ is a Riesz basis for $\h$  satisfying any of the following two conditions:

 $(\mathcal{P}_1)$:\ \ $\forall s\in\mn \ \,  \exists C_s>0 \ : \  |\<g_m, g_n\>|\leq  C_s(1+|m-n|)^{-s}, \ m,n\in\mn$,

 $(\mathcal{P}_2)$:\ \ $\exists s>0 \ \, \exists C_s>0 \ : \    |\<g_m, g_n\>| \leq 
 C_s   \mathrm{e}^{-s |m-n|}, \ m,n\in\mn,
$

\noindent
 then the property $\mprop_{(g_n)}$ is  satisfied.

\vspace{.1in}
Now we apply Theorem \ref{proprb1} to obtain series expansions in the
spaces $\mathcal{S}$ and $\Sigma^\alpha$ (for $\alpha>1/2$), and their duals,  
 via frames localized by the Hermite orthonormal basis and coefficients in the corresponding sequence spaces. 
 Furthermore, we extend the known characterizations of $\mathcal{S}$, $\Sigma^\alpha$, $\alpha>1/2$, and their dual spaces, based on the Hermite basis (see the end of Sec. \ref{sec-Prel}), to characterizations based on a larger class of frame-functions.

\begin{thm} \label{prophb}
Assume that $(e_n)_{n=1}^\infty$  is a sequence with elements from $\mathcal{S}(\mr)$ which is a frame for $L^2(\mr)$ 
and which is  polynomially localized with respect to the Hermite basis $(h_n)_{n=1}^\infty$ with decay $s$ for every $s\in\mn$. Take $(g_n)_{n=1}^\infty:=(h_n)_{n=1}^\infty$. Then $\mprop_{(g_n)}$ and the conclusions in Theorem \ref{proprb1} hold with $X_F$ replaced by $\mathcal{S}$  and $\Theta_F$ replaced by $\bs$. 
 \end{thm}

\begin{thm} \label{prophb2}
 Let $\alpha>1/2$. Assume that a sequence $(e_n)_{n=1}^\infty$ with elements from $\Sigma^\alpha$ is a frame for $L^2(\mr)$ which is exponentially localized with respect to the Hermite basis $\seqgr[h]$.
 Take $(g_n)_{n=1}^\infty:=(h_n)_{n=1}^\infty$.  Then  $\mprop_{(g_n)}$ and  the conclusions in Theorem \ref{proprb1} hold with $X_F$ replaced by $\Sigma^\alpha$  and $\Theta_F$ replaced by $\sspace^{1/(2\alpha)}$. 
 \end{thm}

To illustrate Theorems \ref{prophb} and \ref{prophb2}, below we give an example based on appropriate linear combinations of  Hermite functions.
Further examples are to be given in an extended paper \cite{ps5}.

\begin{ex}  \label{expol}  Let $r\in\mn$ and for $i=1,2,\ldots,r$, take  $\varepsilon_i \geq 0 $  and a sequence $(a_n^i)_{n=1}^\infty$ of complex numbers
 satisfying $|a_n^i|\leq \varepsilon_i$ for $n\geq 2$,  $\sum_{i=1}^r |a^i_1|\leq 1$, and $\sum_{i=1}^r \varepsilon_i <1$. For $n\in\mn$, consider $e_n:=h_n + \sum_{i=1}^r a_n^i h_{n+i},$ which clearly belongs to $\mathcal{S}(\mr)$ and $\Sigma^\alpha$, $\alpha>1/2$. Then the sequence $(e_n)_{n=1}^\infty$ is a Riesz basis for 
$L^2(\mr)$ and it is $s$-localized with respect to the Hermite orthonormal basis $(h_n)_{n=1}^\infty$ for every $s>0$, as well as exponentially localized with respect to $(h_n)_{n=1}^\infty$. 
 \end{ex}

\begin{rem}\label{lag} 
Having in mind the known  expansions of tempered distributions $(\mathcal S(\mr_+))'$ \cite{GTeissier,Duran} 
 and 
Beurling ultradistributions $(G^\alpha_\alpha(\mr_+))'$ \cite{Duran2,JPP}, 
and their test spaces, by the use of the Laguerre orthonormal basis $l_n, n\in\mn,$ and validity of the corresponding properties $\mprop_{(l_n)}$, 
we can transfer the above results to the mentioned classes of distributions and ultradistributions over $\mr_+.$ 
\end{rem}

{\bf Acknowledgements}   
The authors acknowledge support from 
the Multilateral S\&T Danube-Cooperation Project TIFMOFUS (``Time-Frequency Methods for Operators and Function Spaces''; MULT\_DR 01/2017), 
 the Austrian Science Fund (FWF) START-project FLAME ('Frames and Linear Operators for Acoustical Modeling and Parameter Estimation'; Y 551-N13), and the Project 174024 of the Serbian Ministry of Sciences. 
 The second author is grateful for the hospitality of the University of Novi Sad,
  where most of the research on the presented topic was done.

\end{document}